\documentclass{ifacconf}

\usepackage{graphicx}      
\usepackage{natbib}        
\usepackage{algorithmic}

\usepackage{etex}
\usepackage{circuitikz}
\usetikzlibrary{shapes,arrows}
\usetikzlibrary{decorations.pathreplacing}
\usepackage{pgfplots}

\usepackage{epsfig} 
\usepackage{mathptmx} 
\usepackage{times} 
\usepackage{amsmath} 
\usepackage{mathdots}
\usepackage{amssymb}  
\usepackage{float}
\usepackage{placeins}
\usepackage{wrapfig}

\usepackage{mathtools}

\usepackage{epstopdf}
\usepackage{fancyhdr}
\usepackage{booktabs,array}
\usepackage[output-decimal-marker={,}]{siunitx}
\usepackage{color}
\usepackage{empheq}

\newcommand{\al}[1]{\begin{align} #1 \end{align}}
\newcommand{\nn}{\nonumber}
\newcommand{\tr}{\mathrm{tr}}

\newtheorem{theorem}{Theorem}

\newtheorem{remark}{Remark}
\begin{document}
\begin{frontmatter}

\title{Model Predictive Control\\ meets robust Kalman filtering\thanksref{footnoteinfo}} 

\thanks[footnoteinfo]{This work has been partially supported by the FIRB project ``Learning meets time'' (RBFR12M3AC) funded by MIUR.}

\author[First]{Alberto Zenere} 
\author[First]{Mattia Zorzi} 

\address[First]{Dipartimento di Ingegneria dell'€™Informazione, University of Padova, \\via Gradenigo 6/B, 35131 Padova, Italy \\(e-mail: alberto.zenere.1@studenti.unipd.it, zorzimat@dei.unipd.it).}

 \begin{abstract}                

	Model Predictive Control (MPC) is the principal control technique used in industrial applications.
	Although it offers distinguishable qualities that make it ideal for industrial applications, it can be questioned its robustness regarding model uncertainties and external noises.  
	In this paper we propose a robust MPC controller that merges the simplicity in the design of MPC with added robustness. 
	In particular, our control system stems from the idea of adding robustness in the prediction phase of the algorithm through a specific robust Kalman filter recently introduced.
	Notably, the overall result is an algorithm very similar to classic MPC but that also provides the user with the possibility to tune the robustness of the control.
	To test the ability of the controller to deal with errors in modeling, we consider a servomechanism system characterized by nonlinear dynamics.
	
 \end{abstract}

\begin{keyword}
Risk-sensitive filtering,
model predictive control, 
Kalman filtering,
learning control systems,
adaptive control.
\end{keyword}

\end{frontmatter}

\section{Introduction}
In industry the majority of control problems is addressed by using the so called 
	Model Predictive Control, shortly MPC, \cite{Mayne20142967,bib:surveyMPC,camacho2012model}. 
	Indeed the simplicity of the algorithm, based on a straightforward formulation, and the flexibility to deal with constraints, make it the ideal alternative for industrial applications.
	However, it can be questioned its robustness with respect to model uncertainty and noise.
	Although different approaches can be found in literature to overcome this shortcoming (\cite{CALAFIORE2013,bib:RMPCH,bib:surveyMPCrob2,bib:surveyMPCrob,KOTHARE19961361,bib:RMPCimp}), in practice it is usually adopted an {\em ad hoc} MPC tuning \cite[Sec. 1]{MACIEJOWSKI200922}, \cite[Sec. 10]{bib:surveyMPCrob2}. 
	In this respect it is important to develop a control that can be applied to a wide spectrum of applications and that can be tuned even without an in-depth knowledge of the mathematics behind it.

	Therefore, in this paper we present a robust MPC control that allows to add robustness while preserving the intuitive structure of MPC.
    Furthermore we shall see how this algorithm allows the possibility to tune the robustness of the control depending on the needs of the specific application.
    The fundamental observation behind the proposed robust control is that the predictions used in MPC heavily rely on the accuracy of the employed state-space model. 
    Hence the idea is to consider the usual MPC equipped with a robust Kalman filter.  
    More specifically in this paper we explore the use of the robust Kalman filter proposed in \cite{LN2013}, see also \cite{bib:rk,LN2004,HANSEN_SARGENT_2007,ROBUSTNESS_HANSENSARGENT_2008}. 
	Remarkably, the filter admits a Kalman-like structure leading to a simple implementation of the corresponding MPC algorithm and allowing also a reasonably low computational burden. 
	To test the effectiveness of this robust MPC we apply it to control a servomechanism system characterized by nonlinear dynamics and we compare its performances with respect to standard MCP (i.e. MPC equipped with the Kalman filter). Not surprisingly, increasing the robustness of MPC corresponds to increase the energy of the control 
and also to reduce its smoothness.

Finally, we inform the reader that the present paper only reports some preliminary result regarding the robust MPC. In particular, the proof of Theorem \ref{teo1}  and the technical assumptions are omitted and will be published afterwards.

    The outline of this paper is the following.
	In Sec.~\ref{sec:formulation} we briefly review classic MPC formulation, in Sec.~\ref{MPCproposed} we present the main concepts behind the considered robust Kalman filter as well as its applicability combined with MPC, in Sec.~\ref{sec:DC} we introduce the physics of the servomechanism system that we considered and we show different simulations to attest the efficiency of the control system, lastly in Sec.~\ref{sec:conclusions} are drawn the conclusions.

\section{Standard MPC} \label{sec:formulation}

 We consider the discrete-time state-space model
\begin{align}
\Sigma\,:\,\begin{cases}
x_{t+1}& = Ax_{t} + Bu_{t}+Gv_t \\
\quad	y_{t}& = Cx_{t}+Dv_t
\end{cases}
\label{eq:model}
\end{align}
where $x_t\in\mathbb{R}^n$, $u_t\in\mathbb{R}^q$, $y_t\in\mathbb{R}^p$, $v_t\in\mathbb{R}^m$, denote the state, the input, output and unmeasured noise, respectively.
 For simplicity we assume that $GD^T=0$, however the general case can be easily adapted in what follows.
Assume that our task is for the output $y_t$ to follow a certain reference signal $r_t$.
 To this end, we introduce  a quadratic cost function in the form:
\begin{equation}
J_t (\textbf{u}_t,\Sigma)= \sum_{k=1}^{H_p}\Vert \hat{y}_{t+k|t}-r_{t+k}\Vert^2_{Q_k} + \sum_{k=0}^{H_u-1} \Vert \Delta \hat{u}_{t+k|t} \Vert^2_{R_k}
\label{eq:costo}
\end{equation}
where $\hat{y}_{t+k|t}$ represents the prediction of $y_{t+k}$ at time $t$ with $k>0$;  $\Delta\hat{u}_{t+k|t}:=\hat{u}_{t+k|t}-\hat{u}_{t+k-1|t}$ is the predicted variation of the input from time $t+k-1$ to $t+k$ and  $\textbf{u}_t = [u_{t|t}\ .\ .\ .\ u_{t+H_u-1|t}]^T$.
Moreover, $H_p$ denotes the length of the prediction horizon whereas $H_u$ the control one. 
Lastly, $Q_k\in\mathbb{R}^{p\times p}$ and $R_k\in\mathbb{R}^{q\times q}$ indicate the weight matrices for the output prediction errors at time $t$ and the predicted variations of the input at time $t$.
According to the receding horizon strategy, the control input $u_{t|t}$ to apply to $\Sigma$ at time $t$ is extracted
from $\textbf{u}_t$ which is solution to the following open-loop optimization problem  \citep[Chapter 2]{bib:macie}
\al{ \label{optpb} u_{t|t}= \underset{\textbf{u}_t}{\mathrm{argmin}} \,J_t(\textbf{u}_t,\Sigma).}
The well known solution is:

\al{ \label{utt}
u_{t|t} =
[I_q\ 0\ 0\ \cdots] (\Theta^T Q \Theta + R)^{-1}\Theta^T Q(\textbf{r}_t - \Psi \hat x_{t|t})
}
where
\begin{align*}
\Psi =&\left[\begin{array}{cccc} (CA)^T & (CA^2)^T & \ldots & (CA^{H_p})^T \end{array}\right]^T\\
\Theta =&
\begin{bmatrix}
CB & 0_{n\times q} & \cdots & 0_{n\times q}\\
CAB & CB & \cdots & 0_{n\times q}\\
\vdots && \ddots & \vdots  \\
CA^{H_p-1}B && \cdots  & CA^{H_p-H_u-1}B
\end{bmatrix}.\\	
 Q=&\operatorname{diag}(Q_1,\cdots,Q_{H_p})\nn\\
  R=&\operatorname{diag}(R_0,\cdots,R_{H_u-1})\nn\\
\textbf{r}_t =&\left[\begin{array}{cccc} r_{t+1}^T & r_{t+2}^T & \ldots & r_{t+H_p}^T \end{array}\right]^T
\end{align*}
and $\hat x_{t|t}$ denotes the estimate of $x_t$ at time $t$. The latter is typically computed by using the Kalman filter
\al{ \label{KF1}\hat x_{t|t}&= \hat x_{t|t-1} +L_t (y_{t}-C \hat x_{t|t-1})\\
\label{KF2}\hat{x}_{t+1|t} &=  A\hat{x}_{t|t-1} + K_t(y_{t}-C\hat{x}_{t|t-1})+Gu_{t}  \\
\label{KF4bis} L_t &= P_tC^T(CP_tC^T + DD^T)^{-1},\;\; K_t=AL_t\\
\label{KF4}P_{t+1} &= AP_tA^T - K_t(C P_tC^T+DD^T)K_t^T+GG^T  }
where (\ref{KF4bis})-(\ref{KF4}) is the usual Riccati iteration \citep{ferrante2013generalised}. Note that the estimator $\hat x_{t|t}$ is computed by assuming to know the actual underlying model $\Sigma$.

\begin{remark} It is worth noting that the optimization problem (\ref{optpb}) is usually considered with constraints, such as model uncertainty constraints and stability constraints. Such constraints are relatively easy to embed in (\ref{optpb}), on the other hand the price to pay is that the corresponding problem does not admit a closed form solution. Accordingly, its computational burden increases. In what follows we shall continue to consider the unconstrained MPC for the following reasons: first, we will embed the model uncertainty in a different way and second, stability issues are not the focus of the present paper.
\end{remark}
 \section{Robust MPC Proposed}  \label{MPCproposed}
 In practice the actual model, denoted by $\tilde \Sigma$, is unknown and different from the nominal one, denoted by $\Sigma$.
It is then reasonable to assume that we are able to describe this uncertainty, that is we can characterize a set of models $\mathcal{S}$ for which $\tilde \Sigma\in\mathcal{S}$. In the robust MPC formulation the optimization problem (\ref{optpb}) is usually substituted by the mini-max problem \cite[Section 6]{bib:surveyMPCrob2}:
\al{ \label{minimaxpb}u_{t|t}= \underset{\textbf{u}_t}{\mathrm{argmin}} \,\underset{\tilde \Sigma\in\mathcal{S}}{\mathrm{max}} \,J_t(\textbf{u}_t,\tilde \Sigma).}
The latter is sometimes rewritten as a constrained MPC problem.
Many different uncertainty descriptions have been proposed in the literature such as impulse response uncertainty \cite{bib:RMPCimp}, structured feedback uncertainty, \cite{KOTHARE19961361}, polytopic uncertainty, \cite{Angeli20083113}, disturbances uncertainty, \cite{bib:RMPCH}, and probabilistic uncertainty, \cite{LYGEROS_2011}. Finally, in \cite{yang2015risk} it has been proposed a robust MPC wherein the cost function is an exponential-quadratic cost over the state distribution. In this way, large errors are severely penalized.

  The key observation in our MPC formulation is that the standard MPC relies on the the assumption that the actual underlying model is known and thus the Kalman filter (\ref{KF1})-(\ref{KF4}) is designed on it. This assumption, therefore, could deteriorate the performance of MPC when the actual model is different from the nominal one. We propose
  a robust MPC which stems from the idea of building a control system consisting of MPC on one hand, but equipped with a robust state estimator, that takes into account possible differences between the actual and the nominal model, on the other.
In contrast with the usual robust MPC formulation, which is typically based on the mini-max Problem (\ref{minimaxpb}), we consider two independent optimization problems:
\begin{itemize}
\item Robust estimation problem: we want to find a robust estimate $\hat x_{t|t}$ of $x_t$, independently on the the fact that it will be next used to determine the optimal control input $u_{t|t}$. As we will see, this problem is a mini-max problem itself, but its solution gives a robust filter having a Kalman-like structure;
\item Open-loop optimization problem: assuming to have $\hat x_{t|t}$, we want to determine the optimal control input $u_{t|t}$, i.e. it coincides with Problem (\ref{optpb}).
\end{itemize}
It is important to note that this formulation is different from Problem (\ref{minimaxpb}) and represents just an approximation of it. 

In this paper, we consider the robust Kalman filter proposed in \cite{LN2013,bib:rk,LN2004,bib:RKweiner,convtau} which represents the generalization of the risk-sensitive filter \citep{speyer1998,bib:risk,levy2013contraction,bib:RSconvergence}. In the next Section we extend this robust filter to the case wherein an input is present.

  \subsection{Robust Kalman filter with input} \label{sec:Rkalman}
Assume that the nominal model $\Sigma$ is in the form (\ref{eq:model}), where $v_{t}\in\mathbb{R}^m$ is white Gaussian noise (WGN) with E[$v_{t}v^T_{s}$] = $\delta_{t-s}I_m$ ($\delta_t$ represents the Kronecker delta function). Moreover, the noise $v_t$ is independent of the initial state $x_0$, whose nominal distribution is given by $
f_0(x_0)\sim\mathcal{N}(\hat{x}_0,V_0)$. We introduce the random vector $z_t=[\,x_{t+1}^T\,\ y_t^T\,]^T$. At time $t$ the model $\Sigma$ is completely described by the conditional
 probability density of $z_t$ given the measurements $Y_{t-1}:=[\,y_0^T\,y_1^T\,  \ldots \,y_{t-1}^T\,]^T$, denoted by $
{f}_t(z_t|Y_{t-1})$. Note that by construction $ f_t(z_t|Y_{t-1})$  is Gaussian. Let $ \tilde f_t(z_t|Y_{t-1})$ be the conditional probability density of the actual underlying model $\tilde \Sigma$ at time $t$ and assume that it is Gaussian. The discrepancy between $\tilde \Sigma$ and $\Sigma$ at time $t$ can be measured through the Kullback-Leibler divergence \citep{alpha}
\al{\mathcal{D}(\tilde f_t\| f_t)=\int_{\mathbb{R}^{n+p}} \tilde f_t(z_t|Y_{t-1}) \log\left(\frac{\tilde f_t(z_t|Y_{t-1})}{f_t(z_t|Y_{t-1})}\right) \mathrm{d}z_t.}
Thus, we define as the set of all allowable models at time $t$:
\begin{equation}
\mathcal{S}_{t} := \{ \tilde{f}_t(z_t|Y_{t-1})\ |\ \mathcal{D}(\tilde{f}_t || {f}_t) \le c \}
\label{eq:palla}
\end{equation}
that is, at time $t$ the actual model $\tilde \Sigma$ belongs to a ball  about $\Sigma$. The radius of this ball, denoted by $c$, represents the allowable model tolerance and must be fixed {\em a priori}.  Let $\mathcal{G}_t$ denote the class of estimators $g_t$ with finite second-order moments with respect to any probability density $\tilde{f}_{t}(z_{t}|Y_{t-1}) \in \mathcal{S}_{t}$.
\begin{remark} It is worth noting that $\mathcal{S}_t$ can be formed by considering other types of divergence indexes \citep{DUAL,bib:RKdiv,BETA}. \end{remark}

We define as robust estimator of $x_{t+1}$ given $Y_{t}$ the solution to the following minimax problem \begin{equation}
\hat x_{t+1|t} = \underset{g_{t}\in \mathcal{G}_{t}}{\operatorname{argmin}}\ \underset{\tilde{f}_{t}\in {\mathcal{S}}_{t}}{\operatorname{max}}\ \mathbb{E}_{\tilde f_t}[\| x_{t+1}-g_t(y_t)\|^2|Y_{t-1}]
\label{eq:minimax}
\end{equation}
where
\[
\mathbb{E}_{\tilde f_t}[\| x_{t+1}-g_t(y_t)\|^2|Y_{t-1}] = \int_{\mathbb{R}^{n+p}}\ \Vert x_{t+1} - g_t(y_t) \Vert^2 \tilde{f}_t(z_t|Y_{t-1})dz_t
\] is the mean square error of the estimator with respect to the actual model $\tilde \Sigma$. Finally, the robust estimator of $x_t$ given $Y_{t}$, i.e. $\hat x_{t|t}$, is defined as the minimum mean square error estimator based on (\ref{eq:model}) propagating $\tilde f_t(z_t|Y_{t-1})$.

\begin{theorem} \label{teo1}
The robust filter solution to (\ref{eq:minimax}) obeys to the recursion (\ref{KF1})-(\ref{KF2}) where
\al{\label{RK1}L_t&= V_t C^T(CV_t C^T+DD^T)^{-1},\;\; K_t=A L_t  \\
V_t& =(P_t^{-1}-\theta_t I)^{-1}\\
\label{RK1bis}P_{t+1}&= AV_tA^T-K_t(CV_tC^T+DD^T)K_t^T+GG^T }
and $\theta_t>0$ is the unique solution to
\al{\label{RK2}-\log\det(I-\theta_t P_t)^{-1}+\tr[(I-\theta_t P_t)^{-1}-I]=c.}
\end{theorem}

\begin{remark} The robust filter in Theorem \ref{teo1} can be derived by placing the model uncertainty in the transition probability density of $z_t$ given $x_t$, denoted by $\phi_t(z_t|x_t)$. More precisely, the assumption is that the actual transition probability (not necessarily Gaussian) belongs to the ball $\mathcal{T}_t=\{ \tilde \phi_t\, :\, \mathcal{D}(\tilde \phi_t\| \phi_t)\leq c\}$. Such a derivation is more general but less straightforward and is similar to the one in  \cite{LN2013}.
\end{remark}
\begin{remark}
The search of $\theta_t$ satisfying  (\ref{RK2}) can be efficiently performed using the bisection  method.
\end{remark}

\section{MPC of a Servomechanism System}
\label{sec:DC}

\begin{figure*}
\begin{circuitikz}
	
	\draw (0,3) to[V, v_=$V$] (0,0);
	\draw (0,3) to[R, i<^=$I_m$, l=$R$] (3,3);
	\draw (3,3) to[L, l=$L$] (4,3);
	
	\draw (4,3) -- (5,3);
	\draw (5,3) to[V, v_=$E_m$] (5,0);
	\draw (0,0) -- (5,0);
	
	\draw[fill=white] (4.85,0.85) rectangle (5.15,2.15);
	\draw[fill=white] (5,1.5) ellipse (.45 and .45);

	\draw[fill=white] (7,1.5)
	ellipse (.15 and 0.4);
	\draw (7.8,1.5) ellipse (.15 and 0.4);
	\draw[fill=white, color=white] (7.3, 1.1)
	rectangle (7.8, 1.9);
	\draw (7,1.1) -- (7.83,1.1);
	\draw (7,1.9) -- (7.83,1.9);	
	\draw (7.5,1.5) node {$J_m$};
	\draw (7.5,0.7) node {$\beta_m$};
	
	\draw[fill=black] (5.45,1.45) rectangle (7,1.55);
	
	\draw[line width=0.7pt,<-] (5.7,1) arc (-30:30:1);	
	\draw (6.3,1.8) node {$\theta_m$};	
	\draw (6.3,2.3) node {$T_m$};

	\draw[fill=black] (7.95,1.45) rectangle (10.2,1.55);
	
	\draw[line width=0.7pt,->] (8.1,1) arc (-30:30:1);	
	\draw[line width=0.7pt,->] (9.3,1) arc (-30:30:1);	
	\draw (8.6,1.8) node {$T_{f_m}$};	
	\draw (9.7,1.8) node {$T_s$};

	
%
%
%



	\draw[fill=black!50] (10.4,1.49)
	ellipse (.08 and 0.4);
	\draw[fill=black!50, color=black!50] (10.4,1.89)
	rectangle (10.2,1.09);
	\draw[fill=white] (10.2,1.49)
	ellipse (.08 and 0.4);
	\draw (10.2,1.89) -- (10.4,1.89);
	\draw (10.2,1.09) -- (10.4,1.09);

	\draw[fill=black!50] (10.4,0.40)
	ellipse (.13 and 0.67);
	\draw[fill=black!50, color=black!50] (10.4,1.07)
	rectangle (10.2,-0.27);
	\draw[fill=white] (10.2,0.40)
	ellipse (.13 and 0.67);
	\draw (10.2,1.07) -- (10.4,1.07);
	\draw (10.2,-0.27) -- (10.4,-0.27);


	\draw (9.6,0.4) node {$\rho$};
	
	\draw[fill=white] (11.5,0.4)
	ellipse (.15 and 0.4);
	\draw (12.3,0.4) ellipse (.15 and 0.4);
	\draw[fill=white, color=white] (11.8, 0)
	rectangle (12.3, 0.8);
	\draw (11.5,0) -- (12.33,0);
	\draw (11.5,0.8) -- (12.33,0.8);	
	\draw (12,0.4) node {$J_\ell$};
	\draw (12,-0.4) node {$\beta_\ell$};
	
	\draw[fill=black] (10.5,0.35) rectangle (11.5,0.45);
	
	\draw[fill=black] (12.45,0.35) rectangle (13.7,0.45);

	\draw[line width=0.7pt,<-] (12.6,0) arc (-30:30:1);	
	\draw[line width=0.7pt,->] (13.4,0) arc (-30:30:1);	
	\draw (13,0.7) node {$\theta_\ell$};
	\draw (13,1.2) node {$T_d$};		
	\draw (14,0.7) node {$T_{f_\ell}$};
	\end{circuitikz}
\caption{Servomechanism System.}\label{fig:servomotore}
\end{figure*}
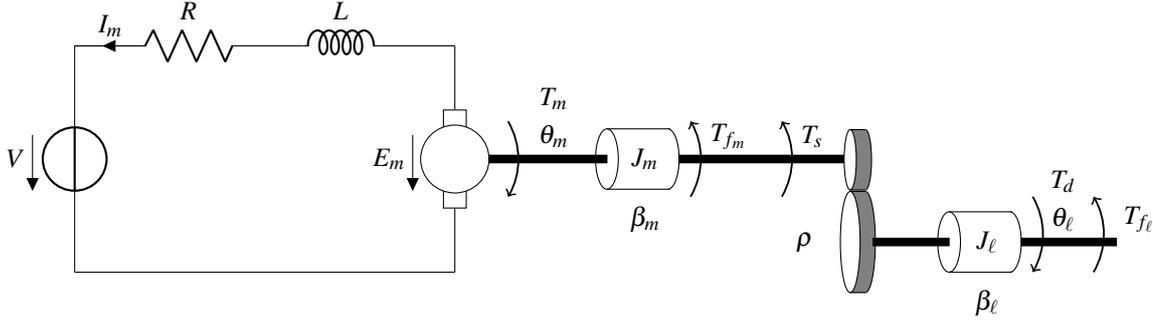

	We consider the servomechanism system presented in \cite{bib:DC}. It consists of a DC-motor, a gear-box, an elastic shaft and a load as depicted in Fig.\ref{fig:servomotore}. The input is the voltage applied $V$ and the output is the load angle $\theta_\ell$.
	Similar mechanisms are often found in the industry for a wide variety of applications.

	The usual approach in MPC is to approximate the system with a linear model, possibly neglecting the nonlinear dynamics. This approximation can be accurate enough as far as conventional control problems are concerned, however it may lead to intolerable errors when the DC motor operates at low speeds and rotates in two directions or when is needed a high precision control. In this situation, indeed, are significant the effects of the Coulomb and the deadzone friction which exhibit nonlinear behaviors in certain regions of operation. Next, we describe the simulation setup and we show the simulation results obtained by applying the standard MPC and the proposed robust MPC.

\subsection{Underlying Nonlinear Model}\label{sec:actual}
\begin{figure}
\begin{center}
\begin{tikzpicture}
\begin{axis} [axis lines=middle,
enlargelimits,
xtick={},
ytick={-9,9},
xticklabels={},
yticklabels={},
xlabel=$\dot{\theta}$,ylabel=$T_f$]
\addplot[domain=0:15, blue, ultra thick] {9+10*exp(-0.2*x)};
\addplot[domain=-15:0, blue, ultra thick] {-9-10*exp(0.2*x)};
\addplot [dotted, domain=0:15] {9};
\addplot [dotted, domain=-15:0] {-9};
\end{axis}
\end{tikzpicture}
\end{center} \caption{Profile of the nonlinear friction model.} \label{fig_prof_nl_attr}
\end{figure}
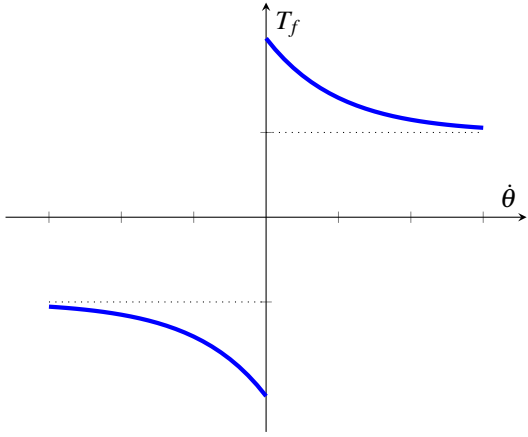

	In this Section we present the complete model including the nonlinear dynamics, that will represent the actual model in our simulations. The equations describing the physics of the system are:
	\begin{align*}
	J_\ell \ddot{\theta}_l &= \rho T_s - \beta_\ell\dot{\theta}_\ell - T_{f_\ell}(\dot{\theta}_\ell) \\
	J_m \ddot{\theta}_m &= T_m - T_s - \beta_m \dot{\theta}_m - T_{f_m}(\dot{\theta}_m) \\
	T_m &= K_t I_m \\
	V &= R I_m + L\dot{I}_m + E_m \\
	E_m &= K_t \dot{\theta}_m \\
	T_s &= \dfrac{k_\theta}{\rho}\left(\dfrac{\theta_m}{\rho}-\theta_\ell\right)
	\end{align*}
where
		$\theta_\ell$ denotes the load angle;
		$\theta_m$ the motor angle;
		$T_m$ the torque generated by the motor;
		$E_m$ the back electromotive force;
		$V$ the motor armature voltage;
		$I_m$ the armature current;
		 $T_s$ the torsional torque;
		$J_\ell$ the load inertia.
Moreover $T_{f_\ell}(\dot{\theta}_\ell)$ and  $T_{f_m}(\dot{\theta}_m)$ represent the Coulomb and the deadzone frictions on the load and on the motor, as described in detail in \cite{bib:attrito}:
\begin{align}
	T_{f_\ell}(\dot{\theta}_\ell) &= \alpha_{\ell 0}\operatorname{sgn}(\dot{\theta}_\ell) + \alpha_{\ell 1}e^{-\alpha_{\ell 2}\vert \dot{\theta}_\ell\vert}\operatorname{sgn}(\dot{\theta}_\ell)
	\label{eq:attrito_l}\\
	T_{f_m}(\dot{\theta}_m) &= \alpha_{m0}\operatorname{sgn}(\dot{\theta}_m) + \alpha_{m1}e^{-\alpha_{m2}\vert \dot{\theta}_m\vert}\operatorname{sgn}(\dot{\theta}_m)
	\label{eq:attrito_m}	
\end{align}
where the function $\operatorname{sgn}$ is defined as:
\[
	\operatorname{sgn}(x) =
	\begin{cases}
		1,\quad x > 0 \\
		0,\quad x = 0 \\
		-1,\quad x < 0.
	\end{cases}
\] The profile of the nonlinear friction model is depicted in Fig. \ref{fig_prof_nl_attr}. Finally, we assume that the motor armature voltage $V$ can take values over the range $\pm 220$ volt. 

  The nominal values of the parameters of the servomechanism system are reported in Table \ref{tab:nom_values}.
  In practice, these values are not accurate because they are difficult to estimate. Accordingly,
 we introduce two kinds of possible parameter perturbations, $\varepsilon_{min}$ and $\varepsilon_{max}$,  expressing the percentage of the relative error that each nominal value can be affected by.
In the actual model every nominal parameter is perturbed in this way: if the nominal value is considered reliable enough it is perturbed by $\pm\varepsilon_{min}=\pm 5\%$, otherwise by $\pm\varepsilon_{max}=\pm 10\%$, see Tab.\ref{tab:real_values}.
\begin{table}[hbt]
	 
\begin{center}
		\caption{Nominal parameters of the servomechanism system.}
	\label{tab:nom_values}
\hfill
	\begin{tabular}{ccl}
		\toprule
		\textit{Symbol} & \textit{Value(MKS)} & \textit{Meaning} \\
		\midrule	
		$L$         &  0       & Armature coil inductance \\
		$J_m$       &    0.5         &   Motor inertia      \\ 		
		$\beta_m$   &    0.1         &   Motor viscous friction coefficient     \\ 		
		$R$         &    20          &   Resistance of armature     \\ 		
		$K_t$       &    10          &   Motor constant      \\ 		
		$\rho$      &    20          &   Gear ratio      \\ 		
		$k_\theta$  &    1280.2      &   Torsional rigidity      \\ 		
		$J_\ell$ &    25    &   Load inertia      \\ 	
		$\beta_\ell$   &    25          &   Load viscous friction coefficient      \\ 												\\						\end{tabular} \hfill
 \end{center}

\end{table}

\begin{table}[hbt]
\caption{Real parameters of the servomechanism system.}
	\label{tab:real_values}
	\makebox[\columnwidth][c]{
	\centering
	\hfill
	\begin{tabular}{ccl}
		\toprule
		\textit{Symbol} & \textit{Value(MKS)} & \textit{Meaning} \\
		\midrule
		$L$ & 0.8 & Armature coil inductance \\	
		$J_m$       &    0.5  (1+$\varepsilon_{max}$)        &   Motor inertia      \\ 		
		$\beta_m$   &    0.1  (1+$\varepsilon_{max}$)        &   Motor viscous friction coefficient     \\ 		
		$R$         &    20   (1+$\varepsilon_{min}$)        &   Resistance of armature      \\ 		
		$K_t$       &    10   (1 + $\varepsilon_{max}$)       &   Motor constant      \\ 		
		$\rho$      &    20   (1+$\varepsilon_{min}$)        &   Gear ratio      \\ 		
		$k_\theta$  &    1280.2 (1+$\varepsilon_{min}$)      &   Torsional rigidity      \\ 		
		${J}_\ell$ &    25 (1-$\varepsilon_{max}$)      &   Load inertia      \\ 	
		$\beta_\ell$   &    25 (1+$\varepsilon_{max}$)          &   Load viscous friction coefficient      \\ 	
		$[\alpha_{\ell_0}\, \alpha_{\ell_1}\, \alpha_{\ell_2}]$   & [0.5\, 10\, 0.5]             &  Load nonlinear friction parameters        \\ 								$[\alpha_{m_0}\, \alpha_{m_1}\, \alpha_{m_2}]$   & [0.1\, 2\, 0.5]             &      Motor nonlinear friction parameters           \\ 	 \\																		
	\end{tabular} \hfill
}
	
\end{table}

\subsection{Linear Model for MPC} \label{secNomimod}
 To obtain a linearized model of the servomechanism system of Section \ref{sec:actual} we eliminate the nonlinear dynamics (\ref{eq:attrito_l})-(\ref{eq:attrito_m}), we set $L=0$ and the saturation on $V$ is removed.
  The dynamic equations resulting from these simplifications are:
\begin{align*}
	\begin{cases}
		J_\ell\ddot{\theta}_\ell &= \rho T_s - \beta_\ell \dot{\theta}_\ell \\
		J_m\ddot{\theta}_m &= T_m - T_s - \beta_m\dot{\theta}_m
	\end{cases}
\end{align*}
wherein we consider the nominal parameters in Tab. \ref{tab:nom_values}.
Defining as state vector  $x_t = [\,\theta_\ell\ \ \dot{\theta}_\ell\ \ \theta_m\ \ \dot{\theta}_m\,]^T$, we obtain a continuous-time state-space linear model of type
\begin{equation*}
\begin{cases}
\dot{x}_t = \bar A x_t + \bar Bu_t + \bar G \dot w_t \\
y_t = \bar Cx_t  + \bar D  \dot w_t \\
\end{cases}
\end{equation*}
where $u_t$ is the motor armature voltage, $y_t$ the load angle, $w_t$ is the normalized Wiener process, and matrices $\bar G$, $\bar D$ are chosen heuristically in such a way to compensate the approximations made before. Lastly this continuous-time model was discretized with sampling time $T=0.1$ s obtaining in this way a discrete-time model of type (\ref{eq:model}).

\subsection{Results}
\label{sec:results}
	In this Section we want to test whether the robust Kalman filter of Section \ref{MPCproposed} does indeed improve the accuracy of MCP. Furthermore it is worth to highlight the effect of the tolerance parameter $c$, introduced in (\ref{eq:palla}), on the overall behavior of the control.	
	
\begin{figure}[hbt]
	\centering
	\includegraphics[width=\columnwidth]{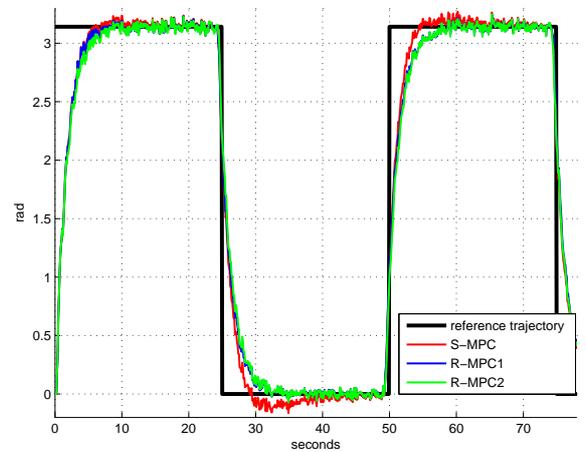}
	\caption{First simulation: load angle when the nominal and the actual model coincide.}
	\label{fig:1_output}
\end{figure}

	To this aim we consider three controllers, they are constituted by MPC coupled with: standard Kalman filter (S-MPC), robust Kalman filter of Section \ref{MPCproposed} with $c=10^{-1}$ (R-MPC1) and $c=1$ (R-MPC2).	
	Thereafter we study their performances in response to the reference trajectory $r_t$, with $t\geq 0$, set as a periodic square wave taking values $0$ and $\pi$ rad with period $T_p=50$ sec and duty cycle$=50$\%, see Fig. \ref{fig:1_output}. The initial state $x_0$ is assumed to be zero mean and with variance equal to the variance of the state-process noise in (\ref{eq:model}). The initial value of the load is $y_0=0$ rad. Regarding the parameters of MPC in the cost function (\ref{eq:costo}) we consider the weight matrices $Q_k = 0.1$, $R_k = 0.1$ and $H_p = 10$, $H_u = 3$.

\begin{figure}[hbt]
	\centering
	\includegraphics[width=\columnwidth]{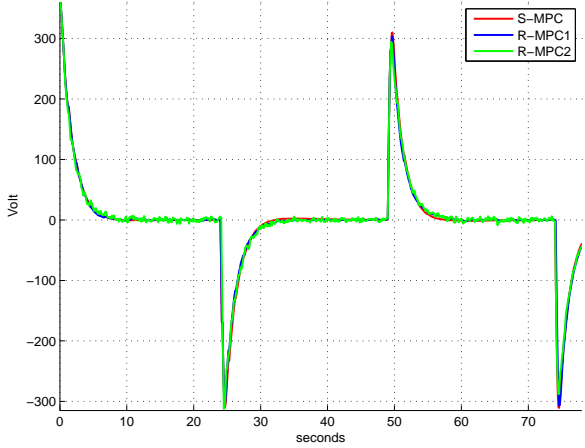}
	\caption{First simulation: voltage applied when the nominal and the actual model coincide.}
	\label{fig:1_input}
\end{figure}

\begin{figure}[hbt]
	\centering
	\includegraphics[width=\columnwidth]{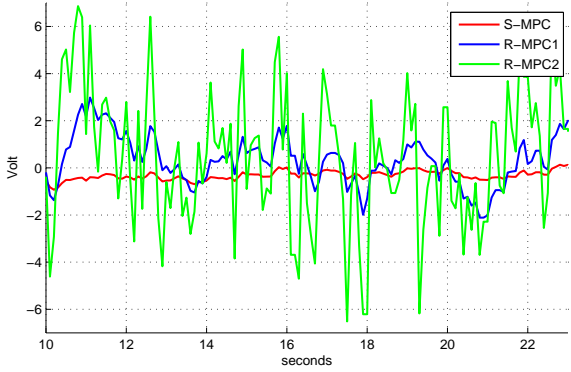}
	\caption{First simulation: voltage applied in the time window [10, 20] sec when the nominal and the actual model coincide.}
	\label{fig:3_output}
\end{figure}

	In the first simulation we consider the ideal situation in which the dynamics of the servomechanism can actually be represented by the linear system computed in Section \ref{secNomimod}, that is the real model coincides with the nominal one.
	As can be seen in Fig.\ref{fig:1_output},
	S-MPC can track the reference trajectory likewise R-MPC1 and R-MPC2. The corresponding voltage applied is depicted in Fig. \ref{fig:1_input};
	in particular in Fig. \ref{fig:3_output} 
	we show the voltage applied during the time window [10, 20] sec, when the load angle already achieved the reference trajectory for all three controllers.
	Here we notice that the input applied by R-MPC1 and R-MPC2 is less smooth than S-MPC and moreover that the smoothness decreases with increasing tolerance. This is not a very surprising fact given that the robust filter is constructed under the idea of considering more uncertainties in the modeling. In particular, the more $c$ is large, the more uncertainty on the nominal model we have.

\begin{figure}[hbt]
	\centering
	\includegraphics[width=\columnwidth]{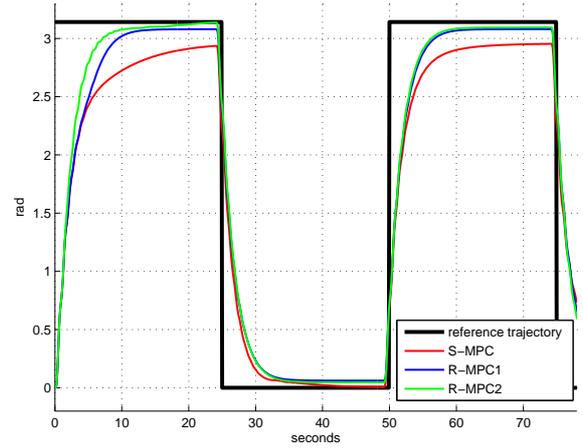}
	\caption{Second simulation: load angle when the actual model differs from the nominal one.}
	\label{fig:2_output}
\end{figure}

\begin{figure}[hbt]
	\centering
	\includegraphics[width=\columnwidth]{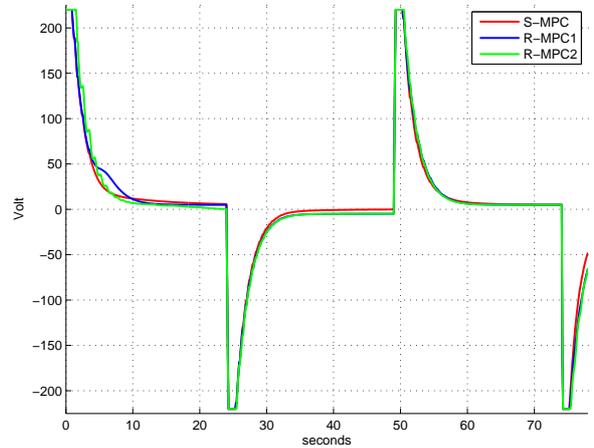}
	\caption{Second simulation: voltage applied when the actual model differs from the nominal one.}
	\label{fig:2_input}
\end{figure}

\begin{figure}[hbt]
	\centering
	\includegraphics[width=0.5\textwidth]{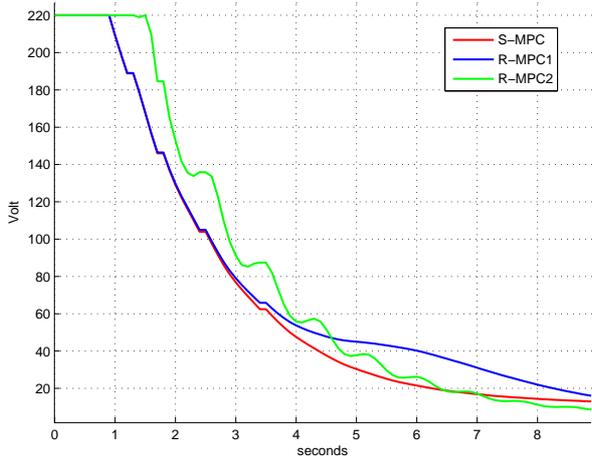}
	\caption{Second simulation: voltage applied in the time window [0,9] sec when the actual model differs from the nominal one.}
	\label{fig:3_input}
\end{figure}

	In the second simulation, the actual model is the one of Section \ref{sec:actual}. Fig. \ref{fig:2_output}
	shows the load angle for S-MPC, R-MPC1 and R-MPC2. 
	It is clear how R-MPC1 and R-MPC2 are able to provide an adequate control whereas S-MPC does not. In particular, R-MPC2 performs slightly better than R-MPC1. In Fig. \ref{fig:2_input} is depicted the corresponding voltage applied 
    and in Fig. \ref{fig:3_input} we show the zoom in the time window [0, 9] sec. 
	It is clear that the price to pay in R-MPC2 is a control with high energy and less smooth.

	To conclude, these simulations highlight how R-MPC provides better performances than S-MPC, under the same fixed weight matrices $Q$ and $R$. More precisely the performance improves as the tolerance parameter $c$ increases, however the drawback is that the control requires more energy and it is less smooth.

 \section{Conclusions}
\label{sec:conclusions}

  Model Predictive Control is the most principled method to address process control. One main challenge that is now presented to the research community is to find a way to gather, in the same controller, the advantages of MPC with robustness properties.
  In this paper we presented an alternative version of MPC based on a robust Kalman filter. Its strong suits comprise an algorithm that is easy to implement and the possibility to tune the robustness of the control.
  To assess the capabilities of this controller, we evaluated its performance with respect to classic MPC. In particular we considered a servomechanism system characterized by nonlinear dynamics. Furthermore, in order to obtain more realistic simulations, a margin of error on the values of the system parameters was introduced.
  Overall, the robust MPC controller proved to be able to consistently compensate errors in modeling more effectively than standard MPC.
  Lastly we evaluated the effect of the tolerance on the performance of the control. As expected, simulations indicate that increasing this parameter allows for improved accuracy on one hand but also requires an input with higher energy.


\begin{thebibliography}{32}
\providecommand{\natexlab}[1]{#1}
\providecommand{\url}[1]{\texttt{#1}}
\providecommand{\urlprefix}{URL }
\expandafter\ifx\csname urlstyle\endcsname\relax
  \providecommand{\doi}[1]{doi:\discretionary{}{}{}#1}\else
  \providecommand{\doi}{doi:\discretionary{}{}{}\begingroup
  \urlstyle{rm}\Url}\fi

\bibitem[{Angeli et~al.(2008)Angeli, Casavola, Franzè, and
  Mosca}]{Angeli20083113}
Angeli, D., Casavola, A., Franzè, G., and Mosca, E. (2008).
\newblock An ellipsoidal off-line {MPC} scheme for uncertain polytopic
  discrete-time systems.
\newblock \emph{Automatica}, 44(12), 3113 -- 3119.

\bibitem[{Banavar and Speyer(1998)}]{speyer1998}
Banavar, R.N. and Speyer, J.L. (1998).
\newblock Properties of risk-sensitive filters/estimators.
\newblock \emph{IEEE Proceedings - Control Theory and Applications}, 145(1),
  106--112.

\bibitem[{Bemporad and Morari(2007)}]{bib:surveyMPCrob2}
Bemporad, A. and Morari, M. (2007).
\newblock Robust model predictive control: a survey.
\newblock \emph{{Lecture Notes in Control and Information Sciences}}, 245,
  207--226.

\bibitem[{Bemporad and Mosca(1998)}]{bib:DC}
Bemporad, A. and Mosca, E. (1998).
\newblock Fulfilling hard constraints in uncertain linear systems by reference
  managing.
\newblock \emph{Automatica}, 34(4), 451--461.

\bibitem[{Boel et~al.(2002)Boel, James, and Petersen}]{bib:risk}
Boel, R.K., James, M.R., and Petersen, I.R. (2002).
\newblock Robustness and risk-sensitive filtering.
\newblock \emph{{IEEE Transactions on Automatic Control}}, 47(3), 451--461.

\bibitem[{Calafiore and Fagiano(2013)}]{CALAFIORE2013}
Calafiore, G.C. and Fagiano, L. (2013).
\newblock Robust model predictive control via scenario optimization.
\newblock \emph{IEEE Transactions on Automatic Control}, 58(1), 219--224.

\bibitem[{Camacho and Bordons(2012)}]{camacho2012model}
Camacho, E.F. and Bordons, C. (2012).
\newblock \emph{Model predictive control in the process industry}.
\newblock Springer Science \& Business Media.

\bibitem[{Chatterjee et~al.(2011)Chatterjee, Hokayem, and
  Lygeros}]{LYGEROS_2011}
Chatterjee, D., Hokayem, P., and Lygeros, J. (2011).
\newblock Stochastic receding horizon control with bounded control inputs: A
  vector space approach.
\newblock \emph{IEEE Transactions on Automatic Control}, 56(11), 2704--2710.

\bibitem[{Ferrante and Ntogramatzidis(2013)}]{ferrante2013generalised}
Ferrante, A. and Ntogramatzidis, L. (2013).
\newblock The generalised discrete algebraic riccati equation in
  linear-quadratic optimal control.
\newblock \emph{Automatica}, 49(2), 471--478.

\bibitem[{Hansen and Sargent(2007)}]{HANSEN_SARGENT_2007}
Hansen, L. and Sargent, T. (2007).
\newblock Recursive robust estimation and control without commitment.
\newblock \emph{Journal of Economic Theory}, 136(1), 1--27.

\bibitem[{Hansen and Sargent(2008)}]{ROBUSTNESS_HANSENSARGENT_2008}
Hansen, L. and Sargent, T. (2008).
\newblock \emph{Robustness}.
\newblock Princeton University Press, Princeton, NJ.

\bibitem[{Jalali and Nadimi(2006)}]{bib:surveyMPCrob}
Jalali, A.A. and Nadimi, V. (2006).
\newblock A survey on robust model predictive control from 1999-2006.
\newblock In \emph{{Internation Conference on Computational Intelligence for
  Modelling Control and Automation and International Conference on Intelligent
  Agents, Web Technologies and Internet Commerce}}.

\bibitem[{Kara and Eker(2004)}]{bib:attrito}
Kara, T. and Eker, I. (2004).
\newblock Nonlinear modeling and identification of a \{DC\} motor for
  bidirectional operation with real time experiments.
\newblock \emph{Energy Conversion and Management}, 45(7-8), 1087 -- 1106.

\bibitem[{Kothare et~al.(1996)Kothare, Balakrishnan, and
  Morari}]{KOTHARE19961361}
Kothare, M.V., Balakrishnan, V., and Morari, M. (1996).
\newblock Robust constrained model predictive control using linear matrix
  inequalities.
\newblock \emph{Automatica}, 32(10), 1361 -- 1379.

\bibitem[{Levy and Nikoukhah(2004)}]{LN2004}
Levy, B.C. and Nikoukhah, R. (2004).
\newblock Robust least-squares estimation with a relative entropy constraint.
\newblock \emph{IEEE Transactions on Information Theory}, 50(1), 89--104.

\bibitem[{Levy and Nikoukhah(2013)}]{LN2013}
Levy, B.C. and Nikoukhah, R. (2013).
\newblock Robust state space filtering under incremental model perturbations
  subject to a relative entropy tolerance.
\newblock \emph{IEEE Transactions on Automatic Control}, 58(3), 682--695.

\bibitem[{Levy and Zorzi(2016)}]{levy2013contraction}
Levy, B.C. and Zorzi, M. (2016).
\newblock A contraction analysis of the convergence of risk-sensitive filters.
\newblock \emph{SIAM Journal on Control and Optimization}, 54(4), 2154--2173.

\bibitem[{Maciejowski(2009)}]{MACIEJOWSKI200922}
Maciejowski, J. (2009).
\newblock Discussion on: âmin-max model predictive control of nonlinear
  systems: A unifying overview on stability".
\newblock \emph{European Journal of Control}, 15(1), 22--25.

\bibitem[{Maciejowski(2001)}]{bib:macie}
Maciejowski, J.M. (2001).
\newblock \emph{Predictive control: with constraints}.
\newblock Pearson education.

\bibitem[{Mayne(2014)}]{Mayne20142967}
Mayne, D.Q. (2014).
\newblock Model predictive control: Recent developments and future promise.
\newblock \emph{Automatica}, 50(12), 2967 -- 2986.

\bibitem[{Nicolao et~al.(1996)Nicolao, Magni, and Scattolini}]{bib:RMPCimp}
Nicolao, G.D., Magni, L., and Scattolini, R. (1996).
\newblock Robust predictive control of systems with uncertain impulse response.
\newblock \emph{Automatica}, 32(10), 1475--1479.

\bibitem[{Orukpe and Jaimoukha(2009)}]{bib:RMPCH}
Orukpe, P.E. and Jaimoukha, I.M. (2009).
\newblock Robust model predictive control based on mixed
  $\mathcal{H}_2/\mathcal{H}_\infty$ control approach.
\newblock \emph{{Proceedings of the European Conference}}, 2223--2228.

\bibitem[{Qin and Badgwell(2003)}]{bib:surveyMPC}
Qin, S.J. and Badgwell, T.A. (2003).
\newblock A survey of industrial model predictive control technology.
\newblock \emph{{Control Engineering Practice II}}, 733--764.

\bibitem[{Yang and Maciejowski(2015)}]{yang2015risk}
Yang, X. and Maciejowski, J. (2015).
\newblock Risk-sensitive model predictive control with gaussian process models.
\newblock \emph{IFAC-PapersOnLine}, 48(28), 374--379.

\bibitem[{Zorzi(2014{\natexlab{a}})}]{BETA}
Zorzi, M. (2014{\natexlab{a}}).
\newblock A new family of high-resolution multivariate spectral estimators.
\newblock \emph{IEEE Transactions on Automatic Control}, 59(4), 892--904.

\bibitem[{Zorzi(2014{\natexlab{b}})}]{alpha}
Zorzi, M. (2014{\natexlab{b}}).
\newblock Rational approximations of spectral densities based on the alpha
  divergence.
\newblock \emph{Mathematics of Control, Signals, and Systems}, 26(2), 259--278.

\bibitem[{Zorzi(2015{\natexlab{a}})}]{DUAL}
Zorzi, M. (2015{\natexlab{a}}).
\newblock An interpretation of the dual problem of the {T}{H}{R}{E}{E}-like
  approaches.
\newblock \emph{Automatica}, 62, 87--92.

\bibitem[{Zorzi(2015{\natexlab{b}})}]{bib:RKdiv}
Zorzi, M. (2015{\natexlab{b}}).
\newblock Multivariate spectral estimation based on the concept of optimal
  prediction.
\newblock \emph{IEEE Transactions on Automatic Control}, 60(6), 1647--1652.

\bibitem[{Zorzi(2015{\natexlab{c}})}]{bib:RKweiner}
Zorzi, M. (2015{\natexlab{c}}).
\newblock On the robustness of the {B}ayes and {W}iener estimators under model
  uncertainty.
\newblock \emph{Submitted, arXiv:1508.01904}.

\bibitem[{Zorzi(2017{\natexlab{a}})}]{convtau}
Zorzi, M. (2017{\natexlab{a}}).
\newblock Convergence analysis of a family of robust {K}alman filters based on
  the contraction principle.
\newblock \emph{Submitted}.

\bibitem[{Zorzi(2017{\natexlab{b}})}]{bib:rk}
Zorzi, M. (2017{\natexlab{b}}).
\newblock Robust {K}alman filtering under model perturbations.
\newblock \emph{IEEE Transactions on Automatic Control}, 62(6).

\bibitem[{Zorzi and Levy(2015)}]{bib:RSconvergence}
Zorzi, M. and Levy, B.C. (2015).
\newblock On the convergence of a risk sensitive like filter.
\newblock In \emph{54th IEEE Conference on Decision and Control (CDC)},
  4990--4995.

\end{thebibliography}
\end{document}